\documentclass[MSNbibl,citesort,seceqn,number,dvips]{arxbj}
\usepackage{upgreek}

\aid{0}
\volume{19}
\issue{5B}
\pubyear{2013}
\firstpage{2250}
\lastpage{2276}
\doi{10.3150/12-BEJ451} 

\makeatletter
\newtheorem{theorem}{Theorem}
\newcommand{\R}{\mathbb{R}}
\def\ve{\varepsilon}
\def\id{\mathrm{id}}
\def\he{\hat\varepsilon}
\def\hell{\hat\ell}
\def\hf{\hat f}
\def\hh{\hat h}
\def\tif{\tilde f}
\def\tiq{\tilde q}
\def\tih{\tilde h}
\def\hr{\hat r}
\def\hD{\hat\Delta}
\def\emp{\mathbb{F}}
\def\hemp{\hat\emp}
\def\hQ{\hat Q}
\def\la{\kappa_a}
\def\1{\mathbf{1}}
\def\ebar{\bar\ve}
\renewcommand{\emptyset}{\varnothing}
\newcommand{\cov}{\operatorname{Cov}}
\def\argmax{\operatorname{\arg\max}\limits}
\def\ahat{\hat a}
\def\E{\mathbb{E}}
\makeatother

\begin{document}
\begin{frontmatter}

\title{Uniform convergence of convolution estimators
for the response density in~nonparametric regression}
\runtitle{Convolution estimators in regression}

\begin{aug}
\author[1]{\fnms{Anton} \snm{Schick}\thanksref{1}\ead[label=e1]{anton@math.binghamton.edu}} \and
\author[2]{\fnms{Wolfgang} \snm{Wefelmeyer}\corref{}\thanksref{2}\ead[label=e2]{wefelm@math.uni-koeln.de}}
\runauthor{A. Schick and W. Wefelmeyer} 
\address[1]{Department of Mathematical Sciences, Binghamton
University, Binghamton, NY 13902, USA.\\\printead{e1}}
\address[2]{Mathematical Institute, University of Cologne, Weyertal
86-90, 50931 K\"oln, Germany.\\\printead{e2}}
\end{aug}

\received{\smonth{12} \syear{2011}}
\revised{\smonth{4} \syear{2012}}

%
\begin{abstract}
We consider a nonparametric regression model $Y=r(X)+\ve$
with a random covariate $X$ that is independent of the error $\ve$.
Then the density of the response $Y$ is a convolution of the densities
of $\ve$ and $r(X)$. It can therefore be estimated by a convolution of kernel
estimators for these two densities, or more generally by a local
von Mises statistic. If the regression function has a nowhere vanishing
derivative, then the convolution estimator converges at a parametric rate.
We show that the convergence holds uniformly,
and that the corresponding process obeys a functional central limit theorem
in the space $C_0(\R)$ of continuous functions vanishing at infinity,
endowed with the sup-norm. The estimator is not efficient.
We construct an additive correction that makes it efficient.
\end{abstract}

%

\end{frontmatter}

\section{Introduction}\label{sec1}

Smooth functionals of densities can be estimated by plug-in estimators,
and densities of functions of two or more random variables can be estimated
by local von Mises statistics.
Such estimators often converge at the parametric rate $n^{1/2}$.
The response density of a nonparametric regression model can be written
in both ways, but it also involves an additional infinite-dimensional parameter,
the regression function. As explained below, this usually leads to
a slower convergence rate of response density estimators,
\emph{except} when the regression function is strictly monotone
in the strong sense that it has a nowhere vanishing derivative.
In the latter case, we can again obtain the rate $n^{1/2}$.

Specifically, consider the nonparametric regression model $Y = r(X)+\ve$
with a one-dimensional random covariate $X$ that is independent
of the unobservable error variable $\ve$.
We impose the following assumptions:

\begin{enumerate}[(G)]
\item[(F)]
The error variable $\ve$ has mean zero,
a moment of order greater than $8/3$,
and a density $f$, and there are bounded and integrable functions
$f'$ and $f''$ such that
$f(z)= \int_{-\infty}^z f'(x)\,\mathrm{d}x$
and $f'(z)= \int_{-\infty}^z f''(x)\,\mathrm{d}x$ for $z\in\R$.
\item[(G)]
The covariate $X$ is \emph{quasi-uniform} on the interval $[0,1]$
in the sense that its density $g$ is bounded and bounded away from zero
on the interval and vanishes outside.
Furthermore, $g$ is of bounded variation.
\item[(R)]
The unknown regression function $r$ is
twice continuously differentiable on $[0,1]$,
and $r'$ is strictly positive on $[0,1]$.
\end{enumerate}
Assume that $(X_1,Y_1),\dots,(X_n,Y_n)$ are $n$ independent copies
of $(X,Y)$.
We are interested in estimating the density $h$ of the response $Y$.
An obvious estimator is the kernel estimator
\[
\frac{1}{n}\sum_{j=1}^nK_b(y-Y_j),
\qquad y\in\R,
\]
where $K_b(t)=K(t/b)/b$ for some kernel $K$ and some bandwidth $b$.
Under the above assumptions on $f$ and $g$, the density $h$ has a
Lipschitz-continuous second derivative as demonstrated in Section \ref{known}.
Thus, if the kernel has compact support and is of order three,
and the bandwidth $b$ is chosen proportional to $n^{-1/7}$,
then the mean squared error of the kernel estimator is of order $n^{-6/7}$.
This means that the estimator has the nonparametric
rate $n^{3/7}$ of convergence.

The above kernel estimator neglects the structure of the regression model.
We shall see that by exploiting this structure one can construct estimators
that have the faster (parametric) rate $n^{1/2}$ of convergence.
For this we observe that the density $h$ is the convolution
of the error density $f$ and the density $q$ of $r(X)$.
The latter density is given by\looseness=-1
\[
q(z) = \frac{g(r^{-1}(z))}{r'(r^{-1}(z))},\qquad z \in\R.
\]\looseness=0
By our assumptions on $r$ and $g$, the density $q$ is quasi-uniform
on the interval $[r(0),r(1)]$, which is the image of $[0,1]$ under $r$.
Furthermore, $q$ is of bounded variation.
The convolution representation $h=f*q$ suggests
a plug-in estimator or \emph{convolution estimator}
$\hh=\hf* \hat q$ based on estimators $\hf$ and $\hat q$ of $f$ and $q$,
for example the kernel estimators
\[
\hf(x) = \frac{1}{n}\sum_{j=1}^n
k_b(x-\he_j)\quad\mbox{and}\quad \hat q(x) =
\frac{1}{n}\sum_{j=1}^n
k_b\bigl(x-\hr(X_j)\bigr),\qquad x\in\R,
\]
with \emph{nonparametric residuals} $\he_j=Y_j-\hr(X_j)$.
Setting $K=k*k$, the convolution estimator $\hh(y)$ has the form
of a \emph{local von Mises statistic}
\[
\frac{1}{n^2} \sum_{i=1}^n\sum
_{j=1}^n K_b\bigl(y-
\he_i-\hr(X_j)\bigr).
\]

In Section \ref{unknown}, we show that the estimator $\hh$
is root-$n$ consistent in the sup-norm and obeys
a functional central limit theorem in the space $C_0(\R)$
of all continuous function on $\R$ that vanish at plus
and minus infinity. As an auxiliary result, Section \ref{known}
treats the case of a \emph{known} regression function $r$.
When $r$ is \emph{unknown}, we estimate it by a local quadratic smoother.
The required properties of this smoother are proved in Section \ref{smooth}.
The convergence rate of $\hh$ follows from a stochastic expansion
which in turn is implied by equations (\ref{e1})--(\ref{e4}).
These equations are proved in Sections \ref{pe1}--\ref{pe4}.

Plug-in estimators in nonparametric settings are often efficient; see,
for example, Bickel and Ritov \cite{BR88}, Laurent \cite{L96}, Chaudhuri \textit
{et al.} \cite{CDS97}
and Efromovich and Samarov \cite{ES00}.
In Section \ref{eff}, we first calculate the asymptotic variance bound
and the efficient influence function for estimators of $h(y)$.
Surprisingly our estimator $\hh(y)$ is not efficient unless
the error distribution happens to be normal. We construct an
additive correction term $\hat C(y)$ such that $\hh(y)-\hat C(y)$
is efficient for $h(y)$. This estimator again obeys a uniform
stochastic expansion and a functional central limit theorem
in $C_0(\R)$. The proof of this result is given in Section \ref{pe5}.

The estimator $\hh$ used here goes back to Frees \cite{F94}.
He observed that densities
of some (known) transformations $T(X_1,\ldots,X_m)$ of $m\geq2$
independent and identically distributed random variables $X_1,\ldots,X_m$
can be estimated pointwise at the parametric rate by a local U-statistic.
Saavedra and Cao \cite{SC00} consider the transformation
$T(X_1,X_2)=X_1+\varphi X_2$ with $\varphi\neq0$.
Schick and Wefelmeyer \cite{SW04b} and \cite{SW07a} obtain this rate
in the sup-norm
and in $L_1$-norms for transformations of the form
$T(X_1,\ldots,X_m)=T_1(X_1)+\cdots+T_m(X_m)$ and $T(X_1,X_2)=X_1+X_2$.
Gin\'e and Mason \cite{GM07} obtain such functional results in $L_p$-norms
for $1\leq p\leq\infty$ and general transformations $T(X_1,\ldots,X_m)$.
The results of Nickl \cite{N07} and \cite{N09} are also applicable in
this context.

The same convergence rates have been obtained for
convolution estimators or local von Mises statistics of the stationary
density of linear processes. Saavedra and Cao \cite{SC99} treat pointwise
convergence for a first-order moving average process.
Schick and Wefelmeyer \cite{SW04a} and \cite{SW04c} consider
higher-order moving average processes and convergence in $L_1$,
and Schick and Wefelmeyer \cite{SW07b} and \cite{SW08a} obtain
parametric rates in the sup-norm and in $L_1$
for estimators of the stationary density of invertible linear processes.
Analogous pointwise convergence results
for response density estimators in nonlinear regression (with
responses missing at random) and in nonparametric regression
are in M\"uller \cite{M12} and St\o ve and Tj\o stheim \cite{ST11},
respectively.
Escanciano and Jacho-Ch\'avez \cite{EJ12} consider the nonparametric regression
model and show uniform convergence, on compact sets, of their
local U-statistic. Their results allow for a multivariate covariate $X$,
but require the density of $r(X)$ to be bounded and Lipschitz.

In the above applications to regression models and time series,
and also in the present paper, the (auto-)regression function
is assumed to have a nonvanishing derivative.
This assumption is essential. Suppose there is a point $x$
at which the regression function behaves like
$r(y)=r(x)+c(y-x)^\nu+\mathrm{o}(|y-x|^\nu)$, for $y$ to the left or
right of $x$,
with $\nu\geq2$. Then the density $q$ of $r(X)$ has a strong peak at $r(x)$.
This slows down the rate of the convolution density estimator
or local von Mises statistic for $h=f*q$. For densities of transformations
$T(X_1,X_2)=|X_1|^\nu+|X_2|^\nu$ of independent and identically
distributed random variables, see Schick and Wefelmeyer \cite{SW08b}
and \cite{SW09} and the review paper by M\"uller \textit{et al.}
\cite{MSW10}.
In their simulations, Escanciano and Jacho-Ch\'avez \cite{EJ12}
consider the
regression function $r(x)=\sin(2\uppi x)$ and a covariate following
a Beta distribution. This choice does not fit their assumptions
because the density of $r(X)$ is neither bounded nor Lipschitz.
Indeed, for $x=1/4$ and $x=3/4$, the regression function
behaves as above with $\nu=2$. In this case, the convolution density
estimator does not have the rate $\sqrt{n}$,
but at best the slower rate $\sqrt{n/\log n}$.

\section{Known regression function}\label{known}

We begin by proving an auxiliary result for the (unrealistic)
case that the regression function $r$ is \emph{known}.
Then we can observe the error $\ve=Y-r(X)$,
and we can apply the results for known transformations
cited in Section \ref{sec1}.
We obtain a root-$n$ consistent estimator of the response density $h$
by the local von Mises statistic
\[
\tih(y)= \frac{1}{n^2} \sum_{i=1}^n
\sum_{j=1}^n K_b\bigl(y-
\ve_i-r(X_j)\bigr),\qquad y\in\R.
\]
In the following, we specify conditions under which the convergence holds
in $C_0(\R)$. We shall assume that $K$ is the convolution $k*k$
for some continuous third-order kernel $k$ with compact support.
Then we can write
\[
\tih(y)= \tif* \tiq(y), \qquad y\in\R,
\]
where
\[
\tif(x)= \frac{1}{n}\sum_{j=1}^nk_b(x-
\ve_j) \quad\mbox{and}\quad \tiq(x)= \frac{1}{n}\sum
_{j=1}^nk_b\bigl(x-r(X_j)
\bigr),\qquad x\in\R.
\]
Setting $f_b=f*k_b$ and $q_b=q*k_b$, we have the decomposition
\[
\tif* \tiq=f_b*q_b + f_b*(
\tiq-q_b)+q_b*(\tif-f_b)+ (\tif
-f_b)*(\tiq-q_b).
\]
Note that $f_b*q_b=f*q*k_b*k_b=h*K_b$.
Since $q$ is of bounded variation and is quasi-uniform on $[r(0),r(1)]$,
we may and do assume that $q$ is of the form
\[
q(x)= \int_{u\le x} \phi(u) \nu(\mathrm{d}u),\qquad x\in\R,
\]
where $\nu$ is a finite measure with $\nu(\R-[r(0),r(1)])=0$,
and $\phi$ is a measurable function such that $|\phi|\le1$.
This allows us to write
\[
h(y) = \int f(y-x)q(x)\,\mathrm{d}x = \int F(y-u) \phi(u) \nu (\mathrm{d}u),
\]
where $F$ is the distribution function corresponding to the
error density $f$. Indeed,
\begin{eqnarray*}
h(y) &=& \int q(y-x) f(x)\,\mathrm{d}x
\\
&=& \int\!\!\int_{u\leq y-x}f(x) \phi(u) \nu(\mathrm{d}u)\,
\mathrm{d}x
\\
&=& \int\!\!\int_{x\leq y-u}f(x)\,\mathrm{d}x \phi(u) \nu (
\mathrm{d}u).
\end{eqnarray*}
The properties of $f$ now yield that $h$ is three times differentiable
with bounded derivatives
%
%
\begin{eqnarray}
\label{h1} h'(y) &=&\int f(y-u) \phi(u) \nu(\mathrm{d}u),\qquad y
\in\R,
\\
\label{h2} h''(y) &=&\int f'(y-u)
\phi(u) \nu(\mathrm{d}u),\qquad y\in\R,
\\
\label{h3} h'''(y)&=&\int
f''(y-u) \phi(u) \nu(\mathrm{d}u),\qquad y\in\R.
\end{eqnarray}
As $k$ is of order three, so is $K$.
Thus, it follows from a standard argument that
\[
\|h*K_b-h\|=\sup_{y\in\R} \bigl|h*K_b(y)-h(y) \bigr| \leq
C b^3
\]
for some constant $C$.

Next, we note that
$f_b *(\tiq-q_b)= H_1*K_b$ with
\[
H_1(y)= \frac{1}{n}\sum_{j=1}^n
\bigl(f\bigl(y-r(X_j)\bigr) - h(y) \bigr),\qquad y\in \R.
\]
Similarly, $q_b*(\tif-f_b)=H_2*K_b$ with
\[
H_2(y)= \frac{1}{n}\sum_{j=1}^n
\bigl(q(y-\ve_j) - h(y) \bigr),\qquad y\in\R.
\]
As shown in Schick and Wefelmeyer \cite{SW07b},
$n^{1/2}H_1$ converges in $C_0(\R)$
to a centered Gaussian process with covariance function
\[
\Gamma_1(s,t)= \cov \bigl(f\bigl(s-r(X)\bigr),f\bigl(t-r(X)\bigr)
\bigr),\qquad s,t\in\R,
\]
and the following approximation holds,
\[
\|H_1*K_b-H_1\| = \mathrm{o}_p
\bigl(n^{-1/2}\bigr).
\]
We can write
\[
H_2(y)= \int \bigl(\emp(y-x)-F(y-x) \bigr) \phi(x) \nu(\mathrm{d}x),
\qquad y\in\R,
\]
where $\emp$ is the empirical distribution function based on the errors
$\ve_1,\dots,\ve_n$,
\[
\emp(t)= \frac{1}{n}\sum_{j=1}^n\1[
\ve_j \leq t], \qquad t\in\R.
\]
Setting $\Delta=n^{1/2}(\emp-F)$ and writing $\|\cdot\|_1$
for the $L_1$-norm, we obtain for each $\delta>0$ the inequalities
\begin{eqnarray*}
T_1(\delta) &=& \sup_{|y_1-y_2|\le\delta} n^{1/2}
\bigl|H_2*K_b(y_1)-H_2*K_b(y_2)
\bigr|
\\[-2pt]
&\le& \sup_{|y_1-y_2|\le\delta} \int\!\!\int \bigl|\Delta(y_1-x-bu)-
\Delta(y_2-x-bu) \bigr| \bigl|K(u)\bigr| \, \mathrm{d}u \nu(\mathrm{d}x)
\\[-2pt]
&\leq&\|K\|_1 \nu(\R) \sup_{|y_1-y_2|\le\delta} \bigl|\Delta(y_1)-
\Delta(y_2)\bigr|.
\end{eqnarray*}
Similarly, we obtain the inequalities
\begin{eqnarray*}
T_2(M)&=&\sup_{|y|>2M} n^{1/2}\bigl|H_2*K_b(y)\bigr|
\\[-2pt]
&\le& \sup_{|y|>2M} \int\!\!\int\bigl|\Delta(y-x-bu)\bigr| \bigl|K(u)\bigr|\, \mathrm{d}u \nu(
\mathrm{d}x)
\\[-2pt]
&\le&\|K\|_1 \nu(\R) \sup_{|y|> M} \bigl|\Delta(y)\bigr|
\end{eqnarray*}
for all $M$ such that $-M<r(0)-bB< r(1)+bB<M$, where the constant $B$
is such
that the interval $[-B,B]$ contains the support of $K$.
From these inequalities, the characterization of compactness
as given in Corollary 4 of Schick and Wefelmeyer \cite{SW07b},
and the properties of the empirical process,
we obtain tightness of the process $n^{1/2}H_2*K_b$ in $C_0(\R)$.
We also have
\[
n^{1/2} \|H_2*K_b-H_2\| \le\|K
\|_1 \nu(\R) \sup_{|y_1-y_2|\le bB} \bigl|\Delta(y_1)-
\Delta(y_2) \bigr|.
\]
It is now easy to conclude that $n^{1/2}H_2*K_b$ converges in
$C_0(\R)$ to a centered Gaussian process with covariance function
\[
\Gamma_2(s,t)= \cov \bigl(q(s-\ve),q(t-\ve) \bigr),\qquad s,t\in\R.
\]


Finally, we have
\[
\bigl\|(\tif-f_b)*(\tiq-q_b)\bigr\| \leq\|\tif-f_b
\|_2 \| \tiq-q_b\|_2 =\mathrm{O}_p
\bigl((nb)^{-1}\bigr),
\]
where $\|\cdot\|_2$ denotes the $L_2$-norm.\eject

The above yield the following result.

%
\begin{theorem}\label{thm.1}
Suppose \emph{(F)}, \emph{(G)} and \emph{(R)} hold, the kernel $K$ is the convolution
$k*k$ of some continuous third-order kernel $k$ with compact support,
and the bandwidth $b$ satisfies $nb^6\to0$ and $nb^2 \to\infty$. Then
$n^{1/2}(\tih- h)$ converges in distribution in the space $C_0(\R)$
to a centered Gaussian process with covariance function $\Gamma_1+\Gamma_2$.
Moreover,
\[
\|\tih- h - H_1-H_2\| = \mathrm{o}_p
\bigl(n^{-1/2}\bigr).
\]
\end{theorem}

\section{Unknown regression function}\label{unknown}

Our main result concerns the case of an \emph{unknown} regression function
$r$. Then we do not observe the random variables $\ve_i=Y_i-r(X_i)$
and $r(X_j)$. In the local von Mises statistic $\tih$ of Section \ref{known},
we therefore replace $r$ by a nonparametric estimator
$\hr$, substitute the residual $\he_i=Y_i-\hr(X_i)$
for the error $\ve_i$, and plug in surrogates $\hr(X_j)$ for $r(X_j)$.
The resulting estimator for $h=f*q$ is then
\[
\hh(y)=\frac{1}{n^2} \sum_{i=1}^n\sum
_{j=1}^n K_b\bigl(y-
\he_i-\hr(X_j)\bigr),\qquad y\in\R.
\]

Our estimator $\hr$ will be a local quadratic smoother.
More precisely, for a fixed $x$ in $[0,1]$, we estimate $r(x)$ by the first
coordinate $\hr(x)=\hat\beta_1(x)$ of the weighted least squares estimator
\[
\hat\beta(x) = \argmax_{\beta} \frac{1}{nc}\sum
_{j=1}^n w\biggl(\frac{X_j-x}{c}\biggr)
\biggl(Y_j-\beta^{\top
}\psi \biggl(\frac{X_j-x}{c} \biggr)
\biggr)^2,
\]
where $\psi(x)=(1,x,x^2)^{\top}$,
the weight function $w$ is a three times continuously
differentiable symmetric density with compact support $[-1,1]$,
and the bandwidth $c$ is proportional to $n^{-1/4}$.
This means that we undersmooth, since an optimal bandwidth for estimating
a twice differentiable regression function is proportional to $n^{-1/5}$.

We assume that $K$ is the convolution $k*k$ for some twice continuously
differentiable third-order kernel $k$ with compact support.
Then we can write our estimator for $h$ as the convolution
\[
\hh(y)= \hf* \hat q(y),\qquad y\in\R,
\]
of the residual-based kernel estimator of $f$,
\[
\hf(x)= \frac{1}{n}\sum_{j=1}^nk_b(x-
\he_j),\qquad x\in\R,
\]
with the surrogate-based kernel estimator of $q$,
\[
\hat q(x)= \frac{1}{n}\sum_{j=1}^nk_b
\bigl(x-\hr(X_j)\bigr),\qquad x\in\R.
\]
Similarly as in Section \ref{known}, we have the decomposition
\[
\hf*\hat q=f_b*q_b + f_b*(\hat
q-q_b)+q_b*(\hf-f_b)+(
\hf_b-f)*(\hat q-q_b).
\]

Let us introduce
\[
\ebar= \frac{1}{n}\sum_{j=1}^n
\ve_j
\]
and
\[
H_3(y)= \frac{1}{n}\sum_{j=1}^n
\ve_j \bigl(f'\bigl(y-r(X_j)
\bigr)-h'(y) \bigr),\qquad y\in\R.
\]
We can write $H_3$ as the convolution $M*f''$
with
\[
M(z)= \frac{1}{n}\sum_{j=1}^n
\ve_j \bigl(\1\bigl[r(X_j)\le z\bigr]-Q(z) \bigr),\qquad
z\in\R,
\]
where $Q$ denotes the distribution function of $r(X)$.
Write $\sigma^2=E[\ve^2]$ for the error variance. Since
$n E[M^2(z)]$ equals $\sigma^2 Q(z)(1-Q(z))$
and $(1-Q)Q$ is integrable, we obtain from Corollary 4
in Schick and Wefelmeyer \cite{SW07b} and the remark after it
that $n^{1/2}H_3$ converges in distribution
in $C_0(\R)$ to a centered Gaussian process with covariance function
$\sigma^2 \Gamma_3$,
where
\[
\Gamma_3(s,t)= \cov \bigl(f'\bigl(s-r(X)
\bigr),f'\bigl(t-r(X)\bigr) \bigr),\qquad s,t\in\R.
\]
Note that $f'$ and $f''$ are bounded and integrable
and therefore square-integrable.

We shall show in Sections \ref{pe1}--\ref{pe4} that
%
%
\begin{eqnarray}
\label{e1} \bigl\|q_b*(\hf-\tif)-\ebar h'\bigr\|&=&
\mathrm{o}_p\bigl(n^{-1/2}\bigr),
\\
\label{e2} \bigl\|f_b*(\hat q-\tiq)+\ebar
h'+J\bigr\| &=& \mathrm{o}_p\bigl(n^{-1/2}\bigr),
\\
\label{e3} \| \hf-f_b\|_2^2
&=& \mathrm{O}_p \biggl(\frac{1}{nb} \biggr),
\\
\label{e4} \| \hat q-q_b
\|_2^2 &=& \mathrm{o}_p(b).
\end{eqnarray}
The last two statements require also $nb^4/\log^4 n \to\infty$.
These four statements and Theorem \ref{thm.1} yield our main result.

%
\begin{theorem}\label{thm.2}
Suppose \emph{(F)}, \emph{(G)} and \emph{(R)} hold, the kernel $K$ is the convolution
$k*k$ of some twice continuously differentiable third-order kernel
$k$ with compact support, and the bandwidth $b$ satisfies
$nb^6\to0$ and $nb^4/\log^4 n\to\infty$.
Let $\hr$ be the local quadratic estimator
for a weight function $w$ that is a three times continuously
differentiable symmetric density with compact support $[-1,1]$,
and for a bandwidth $c$ proportional to $n^{-1/4}$.
Then $n^{1/2}(\hh- h)$ converges in distribution in the space $C_0(\R)$
to a centered Gaussian process with covariance function
$\Gamma=\Gamma_1+\Gamma_2+\sigma^2 \Gamma_3$.
Moreover, we have the uniform stochastic expansion
%
%
\begin{equation}
\label{hexp} \|\hh- h - H_1-H_2-H_3\| =
\mathrm{o}_p\bigl(n^{-1/2}\bigr).
\end{equation}
\end{theorem}

We should point out that
$\Gamma(s,t)= \cov(H(s),H(t))$ for $s,t\in\R$, where
\[
H(y)= f\bigl(y-r(X)\bigr)+q(y-\ve)-\ve \bigl(f'\bigl(y-r(X)
\bigr)-h'(y) \bigr),\qquad y\in\R.
\]

\section{An efficient estimator}\label{eff}

In this section, we treat the question of efficient estimation for $h$.
For the theory of efficient estimation of real-valued functionals
on nonparametric statistical models, we refer to Theorem 2 in Section 3.3
of the monograph by Bickel \textit{et al.} \cite{BKRW98}.
It follows from (\ref{hexp}) that the estimator $\hh(y)$ has
influence function
\[
I_y(X,Y)=q(y-\ve)-h(y)+f\bigl(y-r(X)\bigr)-h(y) - \ve
\bigl(f'\bigl(y-r(X)\bigr)-h'(y) \bigr).
\]
We shall now show that this differs in general from the efficient
influence function. The latter can be calculated as the projection of
$I_{y}(X,Y)$ onto the tangent space of the nonparametric regression model
considered here. The tangent space consists of all functions of the form
\[
\alpha(X)+\beta(\ve)+ \gamma(X)\ell(\ve),
\]
where the function $\alpha$ satisfies
$\int\alpha(x)g(x)\,\mathrm{d}x=0$ and $\int\alpha^2(x)g(x)\,
\mathrm{d}x <\infty$,
the function $\beta$ satisfies $\int\beta(z)f(z)\,\mathrm{d}z=0=\int z\beta
(z)f(z)\,\mathrm{d}y$
and $\int\beta^2(z)f(z)\,\mathrm{d}z<\infty$,
and\vspace*{1pt} the function $\gamma$ satisfies $\int\gamma^2(x) g(x)\,\mathrm
{d}x<\infty$;
see Schick \cite{S93} for details.
The projection of the influence function onto the tangent space is
\begin{eqnarray*}
I^*_y(X,Y) &=& \bigl[f\bigl(y-r(X)\bigr)-h(y)\bigr] + \bigl[q(y-
\ve)-h(y) - d(y) \ell(\ve)\bigr]
\\
&&{} + \biggl[d(y)-\frac{1}{J} \bigl(f'\bigl(y-r(X)
\bigr)-h'(y) \bigr) \biggr] \ell(\ve).
\end{eqnarray*}
%
Here $\ell=-f'/f$ denotes the score function for location,
$J= \int\ell^2(y)f(y)\,\mathrm{d}y$ is the Fisher information,
which needs to be finite for efficiency considerations,
and $d(y)$ is the expectation $E[q(y-\ve) \ve]$. For later use, we set
\[
\lambda(y)= \frac{\ell(y)}{J} -y,\qquad y\in\R.
\]
To see that $I_y^*(X,Y)$ is indeed the projection of the influence function
onto the tangent space, we note that $I_y^*(X,Y)$ belongs to the
tangent space
and that the difference
\[
I_y(X,Y)- I_y^*(X,Y) = \bigl(f'
\bigl(y-r(X)\bigr)-h'(y) \bigr) \lambda(\ve)
\]
is orthogonal to the tangent space. For this, one uses the well-known identities
$E[\ell(\ve)]=0$ and $E[\ve\ell(\ve)]=1$.

We have $I_y(X,Y)=I_y^*(X,Y)$ if and only if $\lambda=0$,
which in turn holds if and only if $f$ is a mean zero normal density.
Consequently, our estimator is efficient for normal errors,
but not for other errors.

In order to see why our estimator for $h(y)$ is not efficient
in general, consider for simplicity the case of known $f$ and $g$.
The efficient influence function is then $-f'(y-r(X)) \ell(\ve)/J$.
Thus, an estimator
$\hat h(y)$ of $h(y)$ is efficient if it satisfies
the stochastic expansion
\[
\hat h(y)= h(y)- \frac{1}{n}\sum_{j=1}^n
\frac{1}{J} f'\bigl(y-r(X_j)\bigr) \ell(
\ve_j) + \mathrm{o}_p\bigl(n^{-1/2}\bigr).
\]
A candidate would be obtained by replacing, in the relevant terms on
the right-hand side, the unknown $r$ by an estimator $\hat r$,
resulting in the estimator
\[
\int f\bigl(y-\hat r(x)\bigr) g(x)\,\mathrm{d}x - \frac{1}{n}\sum
_{j=1}^n\frac{1}{J} f'\bigl(y-
\hat r(X_j)\bigr) \ell \bigl(Y_j-\hat
r(X_j)\bigr).
\]
This shows that a correction term to the plug-in estimator
$\int f(y-\hat r(x)) g(x)\,\mathrm{d}x$ is required for efficiency.

In the general situation, with $f$, $g$ and $r$ unknown,
we must construct a stochastic term $\hat C(y)$ such that
%
%
\begin{equation}
\label{hc} \hat C(y)= \frac{1}{n}\sum_{j=1}^n
\bigl(f'\bigl(y-r(X_j)\bigr)-h'(y)
\bigr) \lambda(\ve_j)+\mathrm{o}_p\bigl(n^{-1/2}
\bigr).
\end{equation}
Then the estimator $\hh(y)-\hat C(y)$ has influence function $I_y^*(X,Y)$,
%
%
\begin{equation}
\label{hh} \hh(y)-\hat C(y) = \frac{1}{n}\sum_{j=1}^nI_y^*(X_j,Y_j)
+ \mathrm {o}_p\bigl(n^{-1/2}\bigr),
\end{equation}
and hence is efficient. We shall construct $\hat C(y)$ such that
(\ref{hc}), and hence (\ref{hh}), hold uniformly in $y$.
This implies a functional central limit theorem in $C_0(\R)$ also for
the improved estimator $\hh-\hat C$. We mention that tightness of
$n^{1/2} C$,
with
\[
C(y)= \frac{1}{n}\sum_{j=1}^n
\lambda(\ve_j) \bigl(f'\bigl(y-r(X_j)
\bigr)-h'(y) \bigr),\qquad y\in\R,
\]
is verified by the same argument as used for $n^{1/2} H_3$.\vadjust{\goodbreak}

To construct the correction term, we use sample splitting.
Let $m$ denote the integer part of $n/2$. Let $\hr_1$ and $\hr_2$
denote the local quadratic smoothers constructed from the observations
$(X_1,Y_1),\dots,(X_m,Y_m)$ or $(X_{m+1},Y_{m+1}),\dots,(X_n,Y_n)$,
both with the same bandwidth $c$ as before.
Define residuals $\he_{i,j} = Y_j -\hr_i(X_j)$ for $i=1,2$ and
$j=1,\dots,n$,
and kernel density estimators
\[
\hf_1(z) = \frac{1}{m}\sum_{j=1}^m
\la(z-\he_{1,j}),\qquad \hf_2(z) = \frac{1}{n-m}\sum
_{j=m+1}^n\la(z- \he_{2,j} )
\]
and
\[
\hf_3(z)= \frac1n \sum_{j=1}^m
\la(z-\he_{2,j}) + \frac1n \sum_{j=m+1}^n
\la(z-\he_{1,j}),
\]
where $\la(x)= \kappa(x/a)/a$ for some bandwidth $a$ and a density
$\kappa$ fulfilling Condition K of Schick \cite{S93}, such as the
logistic kernel.
Then we can estimate $\ell(z)$ by
\[
\hell_1(z)= -\frac{\hf_1'(z)}{a+\hf_1(z)} \quad\mbox{and}\quad
\hell_2(z)= -\frac{\hf_2'(z)}{a+\hf_2(z)},
\]
the Fisher information $J$ by
\[
\hat J = \frac{1}{n} \sum_{j=1}^m
\hell_2^2(\he_{1,j}) + \frac1n \sum
_{j=m+1}^n \hell_1^2(
\he_{2,j}),
\]
and $\lambda(z)$ by
\[
\hat\lambda_i(z)= \frac{\hell_i(z)}{\hat J}- z,\qquad i=1,2.
\]
Finally, we take $\hat C(y)=\hat C_1(y)+ \hat C_2(y)$ with
\[
\hat C_1(y)= \frac{1}{n} \sum_{j=1}^m
\Biggl(\hf'_3\bigl(y-\hr_2(X_j)
\bigr)- \frac{1}{m} \sum_{i=1}^m
\hf'_3\bigl(y-\hr_2(X_i)\bigr)
\Biggr) \hat\lambda_2(\he_{1,j})
\]
and
\[
\hat C_2(y)= \frac{1}{n} \sum_{j=m+1}^n
\Biggl(\hf_3'\bigl(y-\hr_1(X_j)
\bigr)- \frac{1}{n-m} \sum_{i=m+1}^n
\hf'_3\bigl(y-\hr_1(X_i)\bigr)
\Biggr) \hat\lambda_1(\he_{2,j}).
\]
We have the following result, which is proved in Section \ref{pe5}.

%
\begin{theorem}\label{thm.3}
Suppose \emph{(F)}, \emph{(G)} and \emph{(R)} hold, $f$ has finite Fisher information $J$,
and the bandwidth $a$ satisfies $a\to0$ and $a^8n\to\infty$.
Then we have the stochastic expansion $\|\hat C - C\|= \mathrm{o}_p(n^{-1/2})$.
\end{theorem}

Theorems \ref{thm.2} and \ref{thm.3} imply that the improved estimator
$\hh-\hat C$ has the uniform stochastic expansion
\[
\sup_{y\in\R} \Biggl|\hh(y) -\hat C(y) -h(y) - \frac{1}{n}\sum
_{j=1}^nI_y^*(Y_j,X_j)
\Biggr| = \mathrm{o}_p\bigl(n^{-1/2}\bigr)
\]
and is efficient.
As mentioned above, if the errors happen to be normally distributed,
then $\lambda=0$. Therefore, $C=0$ so that $\hat C$ collapses
in the sense that $\|\hat C\|= \mathrm{o}_p(n^{-1/2})$.

\section{Properties of the local quadratic smoother}\label{smooth}

The weighted least squares estimator $\hat\beta(x)$ satisfies
the normal equation
\[
\bar W(x) \hat\beta(x) = \bar V(x)
\]
with
\begin{eqnarray*}
\bar W(x) &=& \frac{1}{nc}\sum_{j=1}^n
w\biggl(\frac{X_j-x}{c}\biggr) \psi \biggl(\frac{X_j-x}{c} \biggr)
\psi^{\top} \biggl(\frac
{X_j-x}{c} \biggr),
\\
\bar V(x) &=& \frac{1}{nc}\sum_{j=1}^n
w\biggl(\frac{X_j-x}{c}\biggr)Y_j \psi \biggl(\frac{X_j-x}{c}
\biggr).
\end{eqnarray*}
Subtracting from both sides of the normal equation the term
$\bar W(x) \beta(x)$ with
\[
\beta(x) = \bigl(r(x),cr'(x),c^2r''(x)/2
\bigr)^{\top},
\]
we arrive at the equality
\[
\bar W(x) \bigl(\hat\beta(x) -\beta(x) \bigr) = A(x) + B(x),
\]
where
\begin{eqnarray*}
A(x) &=& \frac{1}{nc}\sum_{j=1}^n w
\biggl(\frac{X_j-x}{c}\biggr)\ve_j \psi \biggl(\frac{X_j-x}{c}
\biggr),
\\
B(x) &=& \frac{1}{nc}\sum_{j=1}^n w
\biggl(\frac{X_j-x}{c}\biggr)R(X_j,x) \psi \biggl(\frac{X_j-x}{c}
\biggr),
\end{eqnarray*}
and
\begin{eqnarray*}
R(X_j,x) &=& r(X_j)-r(x)-r'(x)
(X_j-x)-\frac{1}{2} r''(x)
(X_j-x)^2
\\
&=& \int_0^1 (X_j-x)^2
\bigl(r''\bigl(x+s(X_j-x)
\bigr)-r''(x) \bigr) (1-s)\,\mathrm{d}s.
\end{eqnarray*}
Since $r''$ is uniformly continuous on $[0,1]$, we see that
\[
\sup_{0\le x \le1} \bigl|B(x)\bigr| = \mathrm{o}_p\bigl(c^2
\bigr)= \mathrm{o}_p\bigl(n^{-1/2}\bigr).
\]
It follows from the proof of Lemma 1 in M\"uller \textit{et al.} \cite
{MSW09} that
\[
\sup_{0\le x \le1} \bigl|A(x)\bigr|^2 = \mathrm{O}_p \biggl(
\frac{\log n}{nc} \biggr)
\]
and
\[
\sup_{0\le x \le1} \bigl|\bar W(x)-W(x)\bigr|^2 = \mathrm{O}_p
\biggl(\frac{\log
n}{nc} \biggr)
\]
with\vspace*{-1pt}
\[
W(x) = E\bigl[\bar W(x)\bigr] = \int g(x+cu) \psi(u)\psi^{\top}(u) w(u)
\,\mathrm{d}u.
\]
Since $g$ is quasi-uniform on $[0,1]$, there is an $\eta$
with $0<\eta<1$ for which
%
%
\begin{equation}
\label{eigen} \eta< \inf_{|v|=1} v^{\top} W(x)v \le
\sup_{|v|=1} v^{\top} W(x)v < \frac{1}{\eta}
\end{equation}
holds for all $x$ in $[0,1]$ and all $c<1/2$.
From this we obtain the expansion
\[
\sup_{0\le x\le1} \bigl|\bar W^{-1}(x)-W^{-1}(x)\bigr|^2
= \mathrm{O}_p \biggl( \frac{\log n}{nc} \biggr),
\]
where $M^{-1}$ denotes a generalized inverse of a matrix $M$ if its
inverse does not exist.
Combining the above, we obtain that\vspace*{-1pt}
%
%
\begin{equation}
\label{r1} \sup_{0\le x\le1} \bigl|\hr(x)-r(x) -D(x) \bigl(A(x)+B(x)\bigr) \bigr| =
\mathrm{O}_p \biggl(\frac{\log n}{nc} \biggr),
\end{equation}
where $D(x)$ is the first row of $W^{-1}(x)$. For later use, we note that
$|D(x)|^2 \le3/\eta^2$ for all $x$ in $[0,1]$ and $c\le1/2$.
We also have
%
%
\begin{equation}
\label{r2} \sup_{0\le x\le1} \bigl|\hr(x)-r(x) -\hat\varrho(x) \bigr| =
\mathrm{o}_p\bigl(n^{-1/2}\bigr),
\end{equation}
where\vspace*{-1pt}
\[
\hat\varrho(x)= D(x)A(x) = \frac{1}{nc}\sum_{j=1}^n
w\biggl(\frac
{X_j-x}{c}\biggr)\ve_j D(x) \psi \biggl(
\frac
{X_j-x}{c} \biggr).
\]
It is easy to check that\vspace*{-1pt}
\begin{eqnarray*}
\int\hat\varrho^2(x) g(x)\,\mathrm{d}x &=& \mathrm{O}_p
\biggl(\frac
{1}{nc} \biggr),
\\[-2pt]
\frac{1}{n}\sum_{j=1}^n\hat
\varrho^2(X_j) &=& \mathrm{O}_p \biggl(
\frac
{1}{nc} \biggr),
\\[-2pt]
\sup_{0\le x\le1} \bigl|\hat\varrho(x)\bigr|^2 &=& \mathrm{O}_p
\biggl(\frac{\log
n}{nc} \biggr).
\end{eqnarray*}
Thus, we obtain
%
\begin{eqnarray}
\label{r3} \frac{1}{n}\sum_{j=1}^n
\bigl(\hr(X_j)-r(X_j) \bigr)^2 &=&
\mathrm{O}_p \biggl(\frac{1}{nc} \biggr),
\\[-2pt]
\label{r5} \int
\bigl(\hr(x)-r(x) \bigr)^4 g(x)\,\mathrm{d}x &=& \mathrm{O}_p
\biggl(\frac
{\log
n}{n^2c^2} \biggr).
\end{eqnarray}

Let $\chi$ be a square-integrable function.
Then the function $\gamma$ defined by
\[
\gamma(t)=\int \bigl(\chi(x-t)-\chi(x) \bigr)^2 \,\mathrm{d}x = \int
\bigl(\chi(x+t)- \chi(x) \bigr)^2 \,\mathrm{d}x,\qquad t\in\R,
\]
is bounded by $4\|\chi\|_2^2$ and satisfies $\gamma(t)\to0$ as $t\to0$.
Using this and the fact that $w$ has support $[-1,1]$, we derive
\begin{eqnarray*}
E \biggl[ \biggl(\int \bigl(\chi(X \pm cu)-\chi(X) \bigr) u^i w(u)
\,\mathrm{d}u \biggr)^2 \biggr] &\le& 
E \biggl[\int \bigl(
\chi(X \pm cu)-\chi(X) \bigr)^2 w(u)\,\mathrm{d}u \biggr]
\\[-2pt]
&\le& \|g\| \int\gamma(cu) w(u)\,\mathrm{d}u\\[-2pt]
& \to&0.
\end{eqnarray*}
Applying this with $\chi=g$, we can conclude
\[
E\bigl[\bigl|W(X)-g(X) \Psi\bigr|^2\bigr] \to0,
\]
where $\Psi= \int\psi(u)\psi^{\top}(u) w(u)\,\mathrm{d}u$.
From this and (\ref{eigen}), we derive that
\[
E\bigl[\bigl|g(X) W^{-1}(X)-\Psi^{-1}\bigr|^2\bigr] \to0.
\]
In particular, with $e=(1,0,0)^{\top}$,\vspace*{1pt}
\[
E\bigl[\bigl| g(X) D(X) - e^{\top} \Psi^{-1} \bigr|^2\bigr]
\to0.
\]
Let us set
\begin{eqnarray*}
t(X) &=&\int g(X-cu)D(X-cu)\psi(u) w(u)\,\mathrm{d}u
\\[-2pt]
&=&\int \bigl(g(X-cu)D(X-cu)-g(X)D(X) \bigr)\psi(u)w(u)\,\mathrm{d}u
\\[-2pt]
&&{} + \bigl(g(X)D(X)-e^{\top}\Psi^{-1} \bigr) \Psi e +1.
\end{eqnarray*}
Then we have\vspace*{1pt}
\begin{eqnarray*}
E\bigl[\bigl(t(X)-1\bigr)^2\bigr] &\le& 6 E \biggl[\int
\bigl|g(X-cu)D(x-cu)-g(X)D(X) \bigr|^2 w(u)\,\mathrm{d}u \biggr]
\\
&&{} +2 E \bigl[ \bigl|g(X)D(X)-e^{\top} \Psi^{-1} \bigr|^2
\bigr] |\Psi e|^2\\
& \to&0,
\end{eqnarray*}
since $|gD|$ is square-integrable. This can be used to show that
\[
\int\hat\varrho(x)g(x)\,\mathrm{d}x = \frac{1}{n}\sum
_{j=1}^n\ve_j t(X_j)=
\ebar+ \mathrm{o}_p\bigl(n^{-1/2}\bigr).
\]
In view of (\ref{r2}), this yields
%
%
\begin{equation}
\label{r6} \int \bigl(\hr(x)-r(x) \bigr)g(x)\,\mathrm{d}x = \ebar+ \mathrm
{o}_p\bigl(n^{-1/2}\bigr).
\end{equation}

\section{\texorpdfstring{Proof of (\protect\ref{e1})}{Proof of (3.1)}}\label{pe1}

Since $q$ is of bounded variation,
we can write $q_b*(\hf-\tif)=\hat H_2*K_b$, where
\[
\hat H_2(y)= \frac{1}{n}\sum_{j=1}^n
\int \bigl(\hemp(y-z)-\emp (y-z) \bigr) \phi(z)\nu(z),\qquad y\in\R,
\]
with $\hemp$ denoting the empirical distribution function
based on the residuals $\he_1,\dots,\he_n$,
\[
\hemp(t)= \frac{1}{n}\sum_{j=1}^n
\1[ \he_j \leq t],\qquad t\in\R.
\]
It was shown in M\"uller \textit{et al.} \cite{MSW09} that
\[
\|\hemp-\emp-\ebar f\| = \mathrm{o}_p\bigl(n^{-1/2}\bigr).
\]
From this and the representation (\ref{h1}) of $h'$,
we immediately derive the expansion
\[
\bigl\|\hat H_2-\ebar h'\bigr\|=\mathrm{o}_p
\bigl(n^{-1/2}\bigr).
\]
This lets us conclude that
\[
\bigl\|q_b*(\hf-\tif)-\ebar h'\bigr\| = \mathrm{o}_p
\bigl(n^{-1/2}\bigr).
\]

\section{\texorpdfstring{Proof of (\protect\ref{e2})}{Proof of (3.2)}}\label{pe2}

Since $f'$ and $f''$ are bounded, a Taylor expansion and the bounds
(\ref{r2}) and (\ref{r3}) yield the uniform expansion
\[
\sup_{y\in\R} \Biggl| \frac{1}{n}\sum_{j=1}^n
\bigl(f\bigl(y-\hr (X_j)\bigr)-f\bigl(y-r(X_j)
\bigr)+f'\bigl(y-r(X_j)\bigr) \hat\varrho(X_j)
\bigr) \Biggr| = \mathrm{o}_p\bigl(n^{-1/2}\bigr).
\]
Now set
\begin{eqnarray*}
S_1(y) &=& \frac{1}{n}\sum_{j=1}^nf'
\bigl(y-r(X_j)\bigr) \hat\varrho(X_j),
\\
S_2(y) &=& \int f'\bigl(y-r(x)\bigr)\hat
\varrho(x) g(x)\,\mathrm{d}x,
\\
S_3(y) &=& \frac{1}{n(n-1)} \sum
_{i\ne j} f'\bigl(y-r(X_j)\bigr)
\ve_i D(X_j) v_c(X_i-X_j),
\\
S &=& \frac{1}{n(n-1)} \sum_{i\ne j}
\ve_i \biggl(D(X_j) v_c(X_i-X_j)-
\int D(x)v_c(X_i-x)g(x)\,\mathrm{d}x \biggr)
\end{eqnarray*}
with
\[
v_c(z)= w(z/c) \psi(z/c)/c.
\]
Then we have
\[
\biggl\|S_1-\frac{n-1}{n} S_3 \biggr\| \le\bigl\|f'\bigr\|
\frac{1}{n^2} \sum_{j=1}^n \bigl|
\ve_j D(X_j) v_c(0) \bigr| =
\mathrm{O}_p \biggl(\frac{1}{nc} \biggr). %
\]
In view of $h'=f''*Q$, we have the identity
\[
S_3(y)-S_2(y)-h'(y)S= \int
f''(z) U(y-z)\,\mathrm{d}z
\]
with
\begin{eqnarray*}
&&U(z)= \frac{1}{n(n-1)} \sum_{i\ne j}
\ve_i \biggl( \bigl(\1\bigl[r(X_j)\le z\bigr]-Q(z) \bigr)
D(X_j) v_c(X_i-X_j)
\\
&& \hspace*{98pt} {}- \int \bigl(\1\bigl[r(x)\le z\bigr]-Q(z) \bigr)
D(x)v_c(X_i-x)g(x)\,\mathrm {d}x \biggr).
\end{eqnarray*}
The terms in the sum have mean zero and are uncorrelated,
with second moments bounded by
$\sigma^2 \1[r(0)\le z \le r(1)] E[|D(X_2)v_c(X_1-X_2)|^2]$.
Thus, we have
\[
n(n-1)\int E\bigl[U^2(z)\bigr] \,\mathrm{d}z \le\sigma^2
\bigl(r(1)-r(0)\bigr) E\bigl[\bigl|D(X_2)v_c(X_1-X_2)\bigr|^2
\bigr] = \mathrm{O}(1/c),
\]
from which we derive
\[
\bigl\|S_3-S_2-h'S\bigr\| \le\bigl\|f''
\bigr\|_2 \|U\|_2 = \mathrm{o}_p
\bigl(n^{-1/2}\bigr).
\]
%
Similarly, one has $cn(n-1)E[S^2] = \mathrm{O}(1)$ and obtains
\[
\bigl\|h'S\bigr\|= \mathrm{o}_p\bigl(n^{-1/2}\bigr).
\]
Next we have $S_2=N*f''$, where
\begin{eqnarray*}
N(z)&=&\frac{1}{n}\sum_{j=1}^n
\ve_j \int\1\bigl[r(x)\le z\bigr] D(x)v_c(X_j-x)g(x)
\,\mathrm{d}x
\\
&=& \frac{1}{n}\sum_{j=1}^n
\ve_j \int\1\bigl[r(X_j-cu)\le z\bigr]
g(X_j-cu) D(X_j-cu) \psi(u)w(u)\,\mathrm{d}u
\\
&=& N_1(z)+N_2(z)+N_3(z)+ Q(z)N
\end{eqnarray*}
with
\begin{eqnarray*}
N_1(z) &=& \int\frac{1}{n}\sum_{j=1}^n
\ve_j \bigl(\1\bigl[r(X_j-cu)\le z\bigr]-\1
\bigl[r(X_j)\le z\bigr] \bigr)
\\
&& \hspace*{33pt} {}\times g(X_j-cu)D(X_j-cu)
\psi(u)w(u)\,\mathrm{d}u,
\\
N_2(z) &=& \frac{1}{n}\sum_{j=1}^n
\ve_j \1\bigl[r(X_j)\le z\bigr],
\\
N_3(z) &=& \frac{1}{n}\sum_{j=1}^n
\ve_j \bigl(t(X_j)-1\bigr) \bigl(\1\bigl[r(X_j)
\le z\bigr]-Q(z) \bigr),
\\
N &=& \frac{1}{n}\sum_{j=1}^n
\ve_j \bigl(t(X_j)-1\bigr).
\end{eqnarray*}
It is easy to check that $N_2*f''=\bar\ve h'+H_3$. Recall the identity
$Q*f''=h'$. Using these identities, we see that
\[
\bigl\|S_2-\bar\ve h' - H_3\bigr\| \le
\bigl\|h'\bigr\||N|+\bigl\|f''\bigr\|_2 \bigl(
\|N_1\|_2+\|N_3\|_2 \bigr).
\]
We show now that the right-hand side is of order $\mathrm{o}_p(n^{-1/2})$.
First, we calculate
\[
nE\bigl[N^2\bigr]=\sigma^2 E\bigl[\bigl(t(X)-1
\bigr)^2\bigr] \to0.
\]
Second, using the abbreviation
$T(u,z)=\1[r(X-cu)\le z]-\1[r(X)\le z]$, we have
\begin{eqnarray*}
n \int E\bigl[N_1^2(z)\bigr]\,\mathrm{d}z &=&
\sigma^2 \int E \biggl[ \biggl(\int T(u,z)g(X-cu)D(X-cu)\psi(u)w(u)\,
\mathrm{d}u \biggr)^2 \biggr]\,\mathrm{d}z
\\
& \le&\sigma^2 \int E \biggl[\int \bigl(T(u,z)g(X-cu)D(X-cu)\psi(u)
\bigr)^2 w(u)\,\mathrm{d}u \biggr]\,\mathrm{d}z
\\
& \le&\sigma^2 \int E \biggl[ \int T^2(u,z)\,\mathrm{d}z
\bigl(g(X-cu)D(X-cu)\psi(u) \bigr)^2 \biggr] w(u)\,\mathrm{d}u
\\
& \le&\sigma^2 \int E \bigl[ \bigl|r(X-cu)-r(X) \bigr| \bigl(g(X-cu)D(X-cu)
\psi(u) \bigr)^2 \bigr]w(u)\,\mathrm{d}u
\\
& \to&0.
\end{eqnarray*}
Third, we derive
\begin{eqnarray*}
n \int E\bigl[N_3^2(z)\bigr]\,\mathrm{d}z &=&
\sigma^2 \int E \bigl[\bigl(t(X)-1\bigr)^2 \bigl(\1
\bigl[r(X)\le z\bigr]-Q(z)\bigr)^2 \bigr]\,\mathrm{d}z
\\
& =& \sigma^2 E \biggl[\bigl(t(X)-1\bigr)^2 \int \bigl(\1
\bigl[r(X)\le z\bigr]-Q(z) \bigr)^2 \,\mathrm{d}z \biggr]
\\
&\le&\sigma^2 \bigl(r(1)-r(0) \bigr) E\bigl[\bigl(t(X)-1
\bigr)^2\bigr]\\
& \to&0.
\end{eqnarray*}
We can now conclude that $\|S_2-\bar\ve h' -H_3\|= \mathrm{o}_p(n^{-1/2})$.

The above relations show that $\|R+\bar\ve h'+H_3\|=\mathrm
{o}_p(n^{-1/2})$, where
\[
R(y)= \frac{1}{n}\sum_{j=1}^n
\bigl(f\bigl(y-\hr(X_j)\bigr)-f\bigl(y-r(X_j)\bigr)
\bigr).
\]
Note that $f_b*(\hat q-\tiq)=R*K_b$. Thus, the desired (\ref{e2})
follows from
the bound
\begin{eqnarray*}
&&\bigl\|f_b*(\hat q-\tiq)+\bar\ve h'+H_3\bigr\|
\\
&&\quad \leq\bigl\|\bigl(R+\bar\ve h'+H_3
\bigr)*K_b\bigr\| + \bigl\|\bigl(\bar\ve h'+H_3
\bigr)*K_b - \bar\ve h'-H_3\bigr\|
\\
&&\quad \leq\bigl\|R+\bar\ve h'+H_3\bigr\|\|K\|_1 +
\bigl\|\bigl(\bar\ve h'+H_3\bigr)*K_b - \bar\ve
h'-H_3\bigr\|
\end{eqnarray*}
and the tightness of $n^{1/2}(\bar\ve h'+H_3)$ in $C_0(\R)$.

\section{\texorpdfstring{Proof of (\protect\ref{e3})}{Proof of (3.3)}}\label{pe3}

Without loss of generality, we assume that $c<1/2$.
Then we have the inequality
%
%
\begin{equation}
\label{dvb} \bigl|D(x) v_c(X-x)\bigr| \le\frac{3}{\eta c} w \biggl(
\frac{X-x}{c} \biggr),\qquad 0\le x \le1.
\end{equation}
Let us set $\ahat= \hr-r$, and, for a subset $C$ of $\{1,\dots,n\}$,
\[
\ahat_{C}(x)= \frac{1}{n}\sum_{j=1}^n
\1[j\notin C] \bigl(\ve_j + R(X_j,x) \bigr) D(x)
v_c(X_j-x).
\]
Note that $\ahat_{\emptyset}(x)=D(x)(A(x)+B(x))$.
For $l=1,\dots,n$ with $l\ne C$ we have
\begin{eqnarray*}
\bigl|\ahat_{C\cup{l}}(x) -\ahat_{C}(x) \bigr| \le\frac{1}{n} \bigl|
\ve_l+R(X_l,x)\bigr| \bigl|D(x)\bigr| \bigl|v_c(X_l-x)\bigr|
\le\frac{3}{\eta} \frac{|\ve_l|
+ c^2 \omega(c)}{nc} w \biggl(\frac{X_l-x}{c} \biggr),
\end{eqnarray*}
where
\[
\omega(c)= \sup\bigl\{ \bigl|r''(x)-r''(y)\bigr|
\dvt x,y\in[0,1], |x-y|\le c \bigr\}.
\]
We abbreviate $\ahat_{\{i\}}$ by $\ahat_i$
and $\ahat_{\{i,j\}}$ by $\ahat_{i,j}$.
The above inequality and (\ref{r1}) yield the rates
%
%
\begin{eqnarray}
\label{a1} \frac{1}{n}\sum_{j=1}^n
\bigl(\ahat(X_j) -\ahat_{j}(X_j)
\bigr)^2 &=& \mathrm{O}_p \biggl(\frac{\log^2 n}{n^2c^2}
\biggr),
\\
\label{a2} \frac{1}{n}\sum_{j=1}^n
\int \bigl(\ahat(x) -\ahat_{j}(x) \bigr)^2 g(x)\,
\mathrm{d}x &=& \mathrm{O}_p \biggl( \frac{\log^2 n}{n^2c^2} \biggr),
\\
\label{a3} E \bigl[ \bigl(\ahat_1(X_1)-
\ahat_{1,2}(X_1) \bigr)^2 \bigr] &=&
\mathrm{O}_p \biggl(\frac{1}{n^2c} \biggr).
\end{eqnarray}
Let us now set
\[
\bar T(z) = \frac{1}{n}\sum_{j=1}^nT_j(z,
\ahat) \quad\mbox{and}\quad \bar T_*(z) = \frac{1}{n}\sum
_{j=1}^nT_j(z,\ahat_j),
\]
where
\[
T_j(z,a)= k_b\bigl(z-\ve_j +
a(X_j)\bigr) - \int\!\!\int k_b\bigl(z-y+a(x)\bigr)
f(y)g(x)\,\mathrm{d}y\, \mathrm{d}x
\]
for a continuous function $a$.
It follows from the properties of $k$ that
%
%
\begin{equation}
\label{kexp} \int \Biggl(\frac{1}{m} \sum_{i=1}^m
\bigl(k_b(x-x_i) - k_b(x-y_i)
\bigr) \Biggr)^2 \,\mathrm{d}x \leq b^{-3}
\bigl\|k'\bigr\|_2^2 \frac{1}{m} \sum
_{i=1}^m (x_i-y_i)^2
\end{equation}
for real numbers $x_1,\dots,x_m$ and $y_1,\dots,y_m$.
This inequality and statements (\ref{a1}) and (\ref{a2})
yield the rate
\[
\int \bigl(\bar T(z)-\bar T_*(z) \bigr)^2 \,\mathrm{d}z =
\mathrm{O}_p \biggl(\frac{\log^2 n}{b^3n^2c^2} \biggr) =
\mathrm{o}_p \biggl(\frac{1}{nb} \biggr).
\]
The last step used the fact that
$nc^2b^2/\log^2n$ is of order $n^{1/2} b^2/\log^2 n$ and tends to infinity.
In addition, we have
\[
n E\bigl[\bar T_*^2(z)\bigr]= E\bigl[T_1^2(z,
\ahat_1)\bigr] + (n-1) E\bigl[T_1(z,\ahat_1)T_2(z,
\ahat_2)\bigr].
\]
Conditioning on $\xi=(\ve_2,X_2,\dots,\ve_n,X_n)$, we see that
\[
E\bigl[T_1(z,\ahat_1)T_2(z,
\ahat_{1,2})\bigr] = E\bigl[ T_2(z,\ahat_{1,2})E
\bigl(T_1(z,\ahat_1)|\xi\bigr)\bigr] = 0.
\]
Similarly one verifies that $E[T_1(z,\ahat_{1,2})T_2(z,\ahat_{2})]$
and $E[T_1(z,\ahat_{1,2})T_2(z,\ahat_{1,2})]$ are zero.
An application of the Cauchy--Schwarz inequality shows that
\[
E\bigl[T_1(z,\ahat_1)T_2(z,
\ahat_2)\bigr] = E \bigl[ \bigl(T_1(z,
\ahat_1)-T_1(z,\ahat_{1,2}) \bigr)
\bigl(T_2(z,\ahat_2)-T_2(z,
\ahat_{1,2}) \bigr) \bigr]
\]
is bounded by $E [ (T_1(z,\ahat_1)-T_1(z,\ahat_{1,2})
)^2 ]$
which in turn is bounded by
\[
E \bigl[ \bigl(k_b\bigl(z-\ve_1-\ahat_1(X_1)
\bigr) -k_b\bigl(z-\ve_1-\ahat_{1,2}(X_1)
\bigr) \bigr)^2 \bigr].
\]
With the help of (\ref{a3}) and (\ref{kexp}), we thus obtain the bound
\[
\int E\bigl[\bar T_*^2(z)\bigr]\,\mathrm{d}z \le\frac{\|k\|_2^2}{nb} +
\frac{(n-1)}{nb^3} \bigl\|k'\bigr\|_2^2 E \bigl[
\bigl(\ahat_1(X_1)-\ahat_{1,2}(X_1)
\bigr)^2 \bigr] 
= \mathrm{O} \biggl(\frac{1}{nb}
\biggr).
\]
It follows that we have the rate $nb \|\bar T\|_2^2 =\mathrm{O}_p(1)$.

Now we set
\[
\hf_*(z)= \int\!\!\int k_b\bigl(z-y+\ahat(x)\bigr)f(y)\,\mathrm{d}y
g(x)\,\mathrm{d}x = \int f_b\bigl(z+\ahat(x)\bigr)g(x)\,\mathrm{d}x.
\]
Since $\hf-\hf_*$ equals $\bar T$, we have
%
%
\begin{equation}
\label{f1} \| \hf-\hf_*\|_2^2 = \mathrm{O}_p
\biggl(\frac{1}{nb} \biggr).
\end{equation}
A Taylor expansion yields the bound
\[
\int \biggl(\hf_*(z)-f_b(z)-f_b'(z)\int
\ahat(x) g(x)\,\mathrm {d}x \biggr)^2\, \mathrm{d}z \le
\bigl\|f_b''\bigr\|_2^2
\int\ahat^4(x) g(x)\,\mathrm{d}x. 
\]
We have
$\|f_b'\|_2=\|f'*k_b\|_2 \le\|f'\|_2\|k_b\|_1
= \|f'\|_2 \|k\|_1$ and $\|f_b''\|_2\le\|f''\|_2 \|k\|_1$.
Using these bounds, (\ref{r5}) and (\ref{r6}), we obtain the rate
%
%
\begin{equation}
\label{f2} \|\hf_*- f_b\|_2^2 =
\mathrm{O}_p \biggl(\frac{1}{n} \biggr).
\end{equation}
The desired result (\ref{e3}) follows from (\ref{f1}) and (\ref{f2}).

\section{\texorpdfstring{Proof of (\protect\ref{e4})}{Proof of (3.4)}}\label{pe4}

We assume again that $c<1/2$ and set
\[
\hat q_*(z) = \int k_b\bigl(z-r(x) - \ahat(x)\bigr) g(x)\,
\mathrm{d}x,\qquad T'(z,a) = \int k_b'
\bigl(z-r(x)\bigr) a(x)g(x)\,\mathrm{d}x.
\]
An argument similar to the one leading to (\ref{f1}) yields
%
%
\begingroup
\abovedisplayskip=7pt
\belowdisplayskip=7pt
\begin{equation}
\label{q1} \|\hat q-\hat q_*\|_2^2 =
\mathrm{O}_p \biggl(\frac{1}{nb} \biggr).
\end{equation}
Note that $\|k_b'\|_2^2 = \mathrm{O}(b^{-3})$ and $\|k_b''\|_2^2 = \mathrm{O}(b^{-5})$.
A Taylor expansion and (\ref{r5}) yield
\[
\int \bigl(\hat q_*(z)-q_b(z)- T'(z,\ahat)
\bigr)^2\, \mathrm{d}z \le\bigl\|k_b''
\bigr\|_2^2 \int\ahat^4(x)g(x)\,\mathrm{d}x =
\mathrm{O}_p \biggl( \frac{\log n}{b^5 n^2c^2} \biggr) =
\mathrm{o}_p \biggl(\frac
{1}{nb^3} \biggr).
\]
In view of (\ref{r2}), we find
\[
\int \bigl(T'(z,\ahat)-T'(z,\hat\varrho)
\bigr)^2\, \mathrm{d}z \le\bigl\|k_b'
\bigr\|_2^2 \int \bigl(\ahat(x)-\hat\varrho(x)
\bigr)^2 g(x)\, \mathrm{d}x = \mathrm{o}_p \biggl(
\frac{1}{nb^3} \biggr).
\]
Finally, we write
\begin{eqnarray*}
T'(z,\hat\varrho) &=& \frac{1}{n}\sum
_{j=1}^n\ve_j \int k_b'
\bigl(z-r(x)\bigr)D(x) v_c(X_j-x) g(x)\,\mathrm{d}x
\\[-2pt]
&=& \frac{1}{n}\sum_{j=1}^n
\ve_j \int \bigl(k_b'\bigl(z-r(x)
\bigr)-k'_b\bigl(z-r(X_j)\bigr) \bigr)
D(x) v_c(X_j-x) g(x)\, \mathrm{d}x
\\[-2pt]
&&{} + \frac{1}{n}\sum_{j=1}^n
\ve_j k_b'\bigl(z-r(X_j)
\bigr) t(X_j).
\end{eqnarray*}
In view of (\ref{dvb}), we have the bound
\[
\int\bigl|D(x)v_c(X-x)\bigr| g(x)\,\mathrm{d}x \le\frac{3}{\eta} \|g\|.
\]
This inequality and an application of the Cauchy--Schwarz
inequality yield the bound
\[
n \int E\bigl[\bigl(T'(z,\hat\varrho)
\bigr)^2\bigr]\,\mathrm{d}z \le 2 \sigma^2 \biggl(
\frac{3 \|g\|}{\eta} E[U] + \bigl\|k_b'\bigr\|_2^2
E\bigl[t^2(X)\bigr] \biggr) 
\]
with
\begin{eqnarray*}
U &=&\int\!\!\int \bigl(k_b'\bigl(z-r(x)
\bigr)-k_b'\bigl(z-r(X)\bigr) \bigr)^2
\bigl|D(x)v_c(X-x)\bigr|g(x)\,\mathrm{d}x \,\mathrm{d}z
\\[-2pt]
&\le&\bigl\|k_b''\bigr\|_2^2
\frac{3}{\eta}\int \bigl(r(X)-r(x) \bigr)^2 \frac{1}{c} w
\biggl(\frac{X-x}{c} \biggr) g(x)\,\mathrm{d}x.
\end{eqnarray*}
In the last step we used (\ref{dvb}) and the analog of (\ref{kexp})
with $k_b'$ in place of $k_b$.
Since $r$ is Lipschitz on $[0,1]$,
we obtain $E[U]= \mathrm{O}(b^{-5}c^2)= \mathrm{o}(b^{-3})$.
The above relations show that
%
%
\begin{equation}
\label{q2} \|\hat q_*-q_b\|_2^2 =
\mathrm{O}_p \biggl(\frac{1}{nb^3} \biggr) =
\mathrm{o}_p(b).
\end{equation}
%
The desired (\ref{e4}) follows from (\ref{q1}) and (\ref{q2}).\vadjust{\goodbreak}
\endgroup

\section{\texorpdfstring{Proof of Theorem \protect\ref{thm.3}}{Proof of Theorem 3}}\label{pe5}

It suffices to show that $n^{1/2}\|\hat C_i-C_i\|=\mathrm{o}_p(1)$ for
$i=1,2$, with
\[
C_1(y)= \frac{1}{n} \sum_{j=1}^m
\bigl(f'\bigl(y-r(X_j)\bigr)-h'(y)
\bigr) \lambda(\ve_j)
\]
and $C_2=C-C_1$. Since the two cases are similar, we prove only the
case $i=1$.

We begin by writing $n^{1/2}C_1= N*f''$ and $n^{1/2}\hat C_1= \hat
N*\hf_3''$
where
\[
N(z)= N(z,\lambda)=\frac{1}{\sqrt{n}}\sum_{j=1}^m
\lambda(\ve_j) \bigl(\1\bigl[r(X_j)\le z\bigr] - Q(z)
\bigr)
\]
and
\[
\hat N(z)= \hat N(z,\hat\lambda_2)= \frac{1}{\sqrt{n}}\sum
_{j=1}^m\hat\lambda_{2}(
\he_{2,j}) \bigl(\1\bigl[\hr_2(X_j)\le z\bigr]-
\hQ(z,\hr_2) \bigr)
\]
with
\[
\hQ(z,\rho)= \frac{1}{m}\sum_{j=1}^m
\1\bigl[\rho(X_j) \le z\bigr].
\]
In view of $E [\int N^2(z)\,\mathrm{d}z ]= E[\lambda^2(\ve)] \int
Q(z)(1-Q(z))\,\mathrm{d}z <\infty$
and the bound
\[
n^{1/2} \|\hat C_1-C_1\| \le\| \hat N-N
\|_2 \bigl\|\hf_3''\bigr\|_2
+ \|N\|_2 \bigl\|\hf_3''-f''
\bigr\|_2
\]
it suffices to show
%
%
\begin{equation}
\label{nn} \|\hat N-N\|_2 =\mathrm{o}_p(1)
\end{equation}
and
%
%
\begin{equation}
\label{ff} \bigl\|\hf_3''-f''
\bigr\|_2= \mathrm{o}_p(1).
\end{equation}

Let us first prove (\ref{ff}).
With $\hD_i=\hr_i-r$, we have $\he_{i,j}=\ve_j-\hD_i(X_j)$
for $i=1,2$ and $j=1,\dots,n$.
Then we can write
\[
\hf_3''(z)-f''(z)=(m/n)
D_1(z)+ \bigl(1-(m/n)\bigr) D_2(z)
\]
with
\begin{eqnarray*}
D_1(z) &=& \frac{1}{m}\sum_{j=1}^m
\bigl(\la''\bigl(z-\ve_j+
\hD_2(X_j)\bigr)-f''(z)
\bigr),\\
D_2(z)&=& \frac{1}{n-m}\sum_{j=m+1}^n
\bigl(\la''\bigl(z-\ve_j +
\hD_1(X_j)\bigr)- f''(z)
\bigr).
\end{eqnarray*}
Let $\E_2$ denote the conditional expectation given
$X_{m+1},Y_{m+1},\dots,X_n,Y_n$. Using the square-integrability of $f''$
and a standard argument, we find that
\begin{eqnarray*}
\E_2 \biggl[\int D_1^2(z)\,\mathrm{d}z
\biggr] &\le& m^{-1}\int\bigl(\la''(z)
\bigr)^2\, \mathrm{d}z
\\
&&{} + \int\!\!\int\!\!\int \bigl(f''\bigl(z-
\hD_2(x)-au\bigr)-f''(z)
\bigr)^2 \kappa(u)\,\mathrm{d}u g(x)\,\mathrm{d}x \,\mathrm{d}z
\\
&=& \mathrm{O}\bigl(m^{-1}a^{-5}\bigr) +
\mathrm{o}_p(1).
\end{eqnarray*}
Thus, $\|D_1\|_2=\mathrm{o}_p(1)$. Similarly, one verifies $\|D_2\|_2=\mathrm{o}_p(1)$,
and we obtain (\ref{ff}).

To prove (\ref{nn}), we set
\[
\bar N(z) = \bar N(z,\hat\lambda_2) = \frac{1}{\sqrt{n}}\sum
_{j=1}^m\int\hat\lambda_2
\bigl(y-\hD_2(X_j)\bigr) f(y)\,\mathrm{d}y \bigl(\1\bigl[
\hr_2(X_j)\le z\bigr]- \hQ(z,\hr_2) \bigr)
\]
and shall verify
\[
\| \hat N - \bar N - N\|_2 = \mathrm{o}_p(1) \quad
\mbox{and}\quad\|\bar N \|_2 =\mathrm{o}_p(1).
\]
We can write
\[
\hat N- \bar N-N = \frac{\hat L-\bar L- L}{\hat J} + \biggl(\frac{1}{\hat J} -
\frac{1}{J} \biggr) L -(\hat M-\bar M -M)
\]
with $\hat L(z)=\hat N(z,\hell_2)$, $\bar L(z)= \bar N(z,\hell_2)$,
$L(z)= N(z,\ell)$, $\hat M(z)= \hat N(z,\id)$, $\bar M(z)=\bar
N(z,\id)$
and $M(z)= N(z,\id)$ where $\id$ denotes the identity map on $\R$.
Now let $\E$ denote the conditional expectation
given $X_1,\dots,X_n,Y_{m+1},\dots,Y_n$.
Then we find
%
%
\begin{equation}
\label{l0} \E\bigl(\|\hat L-\bar L-L\|_2^2\bigr) \le
\frac{1}{n} \sum_{j=1}^m \bigl(2
\Lambda(X_j) R_{1,j} + 2J R_{2,j} \bigr)
\end{equation}
with
\begin{eqnarray*}
\Lambda(x)&=& \int \bigl(\hell_2\bigl(y-\hD_2(x)\bigr)-
\ell(y) \bigr)^2 f(y)\,\mathrm{d}y,
\\
R_{1,j} &=& \int \bigl(\1\bigl[\hr_2(X_j)\le z
\bigr] - \hQ(z,\hr_2) \bigr)^2 \,\mathrm{d}z,
\\
R_{2,j} &=& \int \bigl(\1\bigl[\hr_2(X_j)\le z
\bigr] - \hQ(z,\hr_2) -\1\bigl[r(X_j)\le z\bigr]+Q(z)
\bigr)^2 \,\mathrm{d}z.
\end{eqnarray*}
By the properties of the quadratic smoother, we have
%
%
\begin{equation}
\label{qs} \frac{1}{n}\sum_{j=1}^n
\hD_2^2(X_j)= \mathrm{O}_p
\bigl(n^{-3/4}\bigr) \quad\mbox{and thus}\quad \frac{1}{n}\sum
_{j=1}^n\bigl|\hD_2(X_j)\bigr|=
\mathrm{O}_p\bigl(n^{-3/8}\bigr).
\end{equation}
Several applications of the Cauchy--Schwarz inequality yield the bound
\begin{eqnarray*}
\frac{1}{n} \sum_{j=1}^m
R_{2,j} & \le& \biggl(\frac{3}{n} + \frac{3}{m} \biggr)
\sum_{j=1}^m \int \bigl(\1\bigl[
\hr_2(X_j)\le z\bigr] - \1\bigl[r(X_j)\le z
\bigr] \bigr)^2 \,\mathrm{d}z 
\\
&&{} + 3 \int \bigl(\hQ(z,r)-Q(z) \bigr)^2 \,\mathrm{d}z.
\end{eqnarray*}
Now we use the identity $(\1[u\le z]-\1[v\le z])^2 = \1[u<z\le v]$,
valid for $u\le v$, and (\ref{qs}), to conclude
%
%
\begin{equation}
\label{r2j} \frac{1}{n} \sum_{j=1}^m
R_{2,j} \le\frac{6}{m} \sum_{j=1}^m
\bigl|\hD_2(X_j)\bigr|+\mathrm{O}_p
\bigl(n^{-1/2}\bigr) = \mathrm{O}_p\bigl(n^{-3/8}
\bigr).
\end{equation}
Using the above identity and the uniform consistency of $\hr_2$, we obtain
%
%
\begin{equation}
\label{r1j} \max_{1\le j \le m} R_{1,j} \le\max_{1\le j \le m}
\frac{1}{m} \sum_{i=1}^m \int \bigl(
\1\bigl[\hr_2(X_j)\le z\bigr] -\1\bigl[
\hr_2(X_i)\le z\bigr] \bigr)^2 \,\mathrm{d}z
= \mathrm{O}_p(1).
\end{equation}
By Lemma 10.1 in Schick \cite{S93} there is a constant $c_*$ so that
%
%
\begin{eqnarray}
\label{l1} \frac{1}{n} \sum_{j=1}^m
\int \bigl(\hell_2\bigl(y-\hD_2(X_j)\bigr)-
\hell_2(y) \bigr)^2 f(y)\,\mathrm{d}y &\le&
\frac{c_*}{a^4 n} \sum_{j=1}^m
\hD_2^2(X_j), 
\\
\label{l2} \int \bigl(\hell_2(y)-\ell(y)
\bigr)^2 f(y)\,\mathrm{d}y &\le&\frac{c_*}{a^6m} \sum
_{j=m+1}^n \hD_2^2(X_j)
\nonumber
\\[-8pt]
\\[-8pt]
&&{}+ \mathrm{O}_p \biggl(\frac{1}{a^{6}m} \biggr)+
\mathrm{o}_p(1).
\nonumber
\end{eqnarray}
From (\ref{l0})--(\ref{l2}) and $a^8n \to\infty$,
we obtain $\| \hat L- \bar L - L\|_2 =\mathrm{o}_p(1)$.
A similar argument yields $\| \hat M- \bar M - M\|_2 =\mathrm{o}_p(1)$.
Using (\ref{l1}), (\ref{l2}) and the operator $\E$, we obtain
\[
\frac{1}{n} \sum_{j=1}^m \bigl(
\hell_2(\he_{2,j})-\ell(\ve_j)
\bigr)^2 = \mathrm{o}_p(1).
\]
It is now easy to see that $\hat J$ is a consistent estimator of $J$.
This completes the proof of $\|\hat N- \bar N - N\|_2 =\mathrm{o}_p(1)$.

We are left to verify $\|\bar N\|_2=\mathrm{o}_p(1)$.
Using the definition of $\hQ(z,\hr_2)$, we can write
\[
\bar N(z) = \frac{1}{\sqrt{n}}\sum_{j=1}^m
\bigl(\1\bigl[\hr_2(X_j)\le z\bigr]- \hQ(z,
\hr_2) \bigr) \biggl(\hD_2(X_j) +
\frac{1}{\hat J} \hat\omega(X_j) \biggr),
\]
where
\[
\hat\omega(X_j) = \int \bigl(\hell_2
\bigl(y-\hD_2(X_j)\bigr)-\hell_2(y) \bigr)
f(y)\,\mathrm{d}y = \int\hell_2(y) \bigl(f\bigl(y+
\hD_2(X_j)\bigr)-f(y) \bigr)\,\mathrm{d}y.
\]
A Taylor expansion yields
\[
f\bigl(y+\hD_2(X_j)\bigr)-f(y) - \hD_2(X_j)f'(y)
= \hD_2^2(X_j)\int_0^1
(1-s) f''\bigl(y+s\hD_2(X_j)
\bigr)\,\mathrm{d}s.
\]
Since $\hell_2$ is bounded by $c_*/a$, we obtain
\[
\bigl|\hat\omega(X_j)+\hD_2(X_j) \hat
J_2 \bigr| \le\frac{c_*}{a} \hD_2^2(X_j)
\int\bigl|f''(y)\bigr| \,\mathrm{d}y
\]
with $\hat J_2= \int\hell_2(y)\ell(y)f(y)\,\mathrm{d}y= J+\mathrm{o}_p(1)$.
Now set
\begin{eqnarray*}
\hat\Upsilon(z) &=& \frac{1}{\sqrt{n}}\sum_{j=1}^m
\bigl(\1\bigl[\hr_2(X_j)\le z\bigr]- \hQ(z,
\hr_2) \bigr)\hD_2(X_j),
\\
\Upsilon(z) &=& \frac{1}{\sqrt{n}}\sum_{j=1}^m
\bigl(\1\bigl[r(X_j)\le z\bigr]- Q(z) \bigr)\hD_2(X_j).
\end{eqnarray*}
Using the Minkowski inequality and the statements (\ref{qs})--(\ref{r1j}),
we derive
\begin{eqnarray*}
\bigl\|\bar N- (1-\hat J_2/\hat J) \hat\Upsilon\bigr\|_2 &\le&
\frac{c_*\|f''\|_1}{a \sqrt{n}} \sum_{j=1}^m
R_{1,j}^{1/2} \hD_2^2(X_j)
= \mathrm{O}_p\bigl(a^{-1} n^{-1/4}\bigr)=
\mathrm{o}_p(1),
\\
\|\hat\Upsilon- \Upsilon\|_2 &\le& \frac{1}{\sqrt{n}}\sum
_{j=1}^mR_{2,j}^{1/2} \bigl|
\hD_2(X_j)\bigr| \le n^{1/2} \Biggl(
\frac{1}{n} \sum_{j=1}^m
R_{2,j} \frac{1}{n}\sum_{j=1}^m
\hD_2^2(X_j) \Biggr)^{1/2} =
\mathrm{o}_p(1).
\end{eqnarray*}
Using the inequality $|\1[r(x)\le z]-Q(z)|\le\1[r(0)\le z\le r(1)]$,
valid for all $0\le x\le1$ and $z\in\R$, we obtain
\begin{eqnarray*}
\E_2\bigl[\bigl \|\Upsilon-\E_2[\Upsilon]
\bigr\|_2^2\bigr] &\le&\frac{m}{n} \int\!\!\int \bigl(\1
\bigl[r(x)\le z\bigr]-Q(z) \bigr)^2 \hD_2^2(x)
g(x)\,\mathrm{d}x \,\mathrm{d}z
\\
&\le& \bigl(r(1)-r(0) \bigr) \int\hD_2^2(x)g(x)\,
\mathrm{d}x =\mathrm{o}_p(1).
\end{eqnarray*}
Now introduce
\[
I(z,\rho)= \int \bigl(\1\bigl[r(x)\le z\bigr] -Q(z) \bigr) \rho(x)g(x)\,
\mathrm{d}x.
\]
Then we have $\E_2[\Upsilon(z)]= n^{-1/2} m I(z,\hD_2)$.
In view of the above and $1-\hat J_2/\hat J=\mathrm{o}_p(1)$,
the desired property $\|\bar N\|_2=\mathrm{o}_p(1)$ will follow if we show
$\| I(\cdot,\hD_2)\|_2 = \mathrm{O}_p(n^{-1/2})$.
The latter is equivalent to showing
$\|I(\cdot,\hr-r)\|_2 = \mathrm{O}_p(n^{-1/2})$. In view of (\ref{r2}),
we have
\[
\bigl\| I(\cdot, \hr-r) - I(\cdot,\hat\varrho)\bigr\|_2 = \bigl\|I(\cdot,\hr -r-
\hat\varrho)\bigr\|_2 =\mathrm{o}_p\bigl(n^{-1/2}
\bigr).
\]
We can express $I(z,\hat\varrho)$ as the average
\[
\frac{1}{n}\sum_{j=1}^n
\ve_j \tau(z,X_j)
\]
with
\begin{eqnarray*}
\tau(z,X_j) &=& \int \bigl(\1\bigl[r(x)\le z\bigr]-Q(z) \bigr)
\frac{1}{c} w \biggl(\frac
{X_j-x}{c} \biggr) D(x) \psi \biggl(
\frac{X_j-x}{c} \biggr) g(x)\,\mathrm{d}x
\\
&=& \int \bigl(\1\bigl[r(X_j-cu)\le z\bigr]-Q(z) \bigr) w(u)
D(X_j-cu)\psi(u) g(X_j-cu)\,\mathrm{d}u.
\end{eqnarray*}
Since $|\tau(z,X_j)|$ is bounded by a constant times $\1[r(0)\le z \le r(1)]$,
we conclude
\[
n E\bigl[\bigl\| I(\cdot, \hat\varrho)\bigr\|_2^2\bigr] = \int
\sigma^2 E\bigl[\tau^2(z,X)\bigr] \,\mathrm{d}z =
\mathrm{O}(1).
\]
The above shows that $\|I(\cdot,\hr-r)\|_2 = \mathrm{O}_p(n^{-1/2})$,
and the proof is finished.

\section*{Acknowledgements}

We thank the referee for suggesting to discuss the question of efficiency.
This resulted in adding the present Sections \ref{eff} and \ref{pe5}.
The research of Anton Schick was supported in part by
NSF Grant DMS 09-06551.


%

\printhistory

\end{document}